%% file: arxiv_conf_GRAPHEM.tex
\documentclass{article}
\usepackage{spconf,amsmath,graphicx,bm}
\usepackage{amssymb}
\usepackage{mathtools} 
\usepackage[bottom]{footmisc}
\usepackage{breqn}
\usepackage{threeparttable}
\usepackage[square,numbers,sort&compress]{natbib}
\usepackage{url}
\usepackage{algorithm}
\usepackage{multirow}
\usepackage[normalem]{ulem}

\def\x{{\mathbf x}}
\def\y{{\mathbf y}}
\def\q{{\mathbf q}}
\def\r{{\mathbf r}}
\def\A{{\mathbf A}}
\def\B{{\mathbf B}}
\def\C{{\mathbf C}}
\def\H{{\mathbf H}}
\def\R{{\mathbf R}}
\def\Q{{\mathbf Q}}
\def\Y{{\mathbf Y}}
\def\V{{\mathbf V}}
\def\G{{\mathbf G}}
\def\X{{\mathbf X}}
\def\Z{{\mathbf Z}}
\def\m{{\mathbf m}}
\def\z{{\mathbf z}}
\def\P{{\mathbf P}}
\def\S{{\mathbf S}}
\def\Id{{\textbf{Id}}}

\newcommand{\Sigmab}{{\bm \Sigma}}
\newcommand{\Phib}{{\bm \Phi}}
\newcommand{\thetab}{{\bm \theta}}
\usepackage{color}

\title{GraphEM: EM algorithm for blind Kalman filtering under graphical sparsity constraints}
\name{\'Emilie Chouzenoux$^{1}$ and V\'ictor Elvira$^{2}$}

\address{$^{1}$ Universit\'e Paris-Saclay, CentraleSup\'elec, Inria, CVN, Gif-sur-Yvette, France \\
    $^{2}$ School of Mathematics, University of Edinburgh, United Kingdom
    \thanks{V.E. and \'E.C. acknowledge support from the Agence Nationale de la
Recherche of France under PISCES (ANR-17-CE40-0031-01) and MAJIC
(ANR-17-CE40-0004-01) projects.}}
    
\begin{document} 
\maketitle
\begin{abstract} 
Modeling and inference with multivariate sequences is central in a number of signal processing applications such as acoustics, social network analysis, biomedical, and finance, to name a few. The linear-Gaussian state-space model is a common way to describe a time series through the evolution of a hidden state, with the advantage of presenting a simple inference procedure due to the celebrated Kalman filter. A fundamental question when analyzing multivariate sequences is the search for relationships between their entries (or the modeled hidden states), especially when the inherent structure is a non-fully connected graph. In such context, graphical modeling combined with parsimony constraints allows to limit the proliferation of parameters and enables a compact data representation which is easier to interpret by the experts. In this work, we propose a novel expectation-minimization algorithm for estimating the linear matrix operator in the state equation of a linear-Gaussian state-space model. Lasso regularization is included in the M-step, that we solved using a proximal splitting Douglas-Rachford algorithm. Numerical experiments illustrate the benefits of the proposed model and inference technique, named GraphEM, over competitors relying on Granger causality. 
\end{abstract}
\section{Introduction}
\label{sec:intro} 
\input{00_introduction}

\section{BACKGROUND}
\label{sec:format}
\input{01_background}

\section{PROPOSED GRAPHEM APPROACH}
\label{sec:proposed}
\input{02_graphem}

\section{EXPERIMENTAL EVALUATION}
\label{sec:experiment}
\input{03_experiments}

\section{CONCLUSION}
\label{sec:conclusion}
This work proposes the GraphEM algorithm, an expectation-minimization method for the estimation of the linear operator that encodes the hidden state relationships in a linear-Gaussian state-space model. Due to the use of a Lasso penalization term, we incorporate a sparsity constraint in our model. This is particularly well suited for modeling and representing the state entries interactions as a compact and interpretable graph. A proximal splitting algorithm is employed for solving the inner minimization problem of our EM approach. Numerical results illustrate the great performance of the method, when compared to several techniques for graphical modeling.  


\bibliographystyle{IEEEbib}
\footnotesize

\end{document}

%% file: 00_introduction.tex
State-space models allow for the statistical description of systems in countless applications in science and engineering. They model a hidden state that evolves over (discrete) time, which is only partially observed through transformed and noisy partial observations. In its simple form, the state-space models have a Markovian dependency in the hidden-state, and each observation depends on the current state, being conditionally independent on the rest of states.

In this framework, a common statistical task is the estimation of the hidden state conditioned to the available observations. When the estimation is done in a probabilistic manner, this is commonly known as Bayesian filtering. Similarly, when the estimation of each state is done using the whole sequence of observations, the task is called Bayesian smoothing. The linear-Gaussian state-space model has been widely used in a plethora of problems and it admits an exact computation of both filtering and smoothing distribution through the Kalman filter and the Rauch-Tung-Striebel (RTS) smoother \cite[Chapter 8]{Sarkka}. For non-linear state-space models, there exist many other methods such as the extended KF (EKF), the unscented KF (UKF) \cite{Julier04}, or the particle filters (PF) \cite{Djuric04}. However, all these algorithms need the model parameters to be known. In the Kalman filter literature, there exist methods for the estimation of the static parameters, based on the expectation-minimization (EM) algorithm or optimization-based methods \cite[Chapter 12]{Sarkka}.

In this paper, we propose a method called GraphEM for estimating the linear matrix of the state equation of a linear-Gaussian state-space model. While this parameter is arguably the hardest to be estimated, it provides very valuable information about the hidden process, not only for inference purposes but also for understanding the uncovered relations among the state dimensions, in the line of graphical modeling methods for time series~\cite{Eichler2012,Bach04,Barber10}. Such representation of multivariate sequences interactions has applications in several domains such as biology \cite{Pirayre,Luengo19}, social network analysis \cite{Ravazzi} and neurosciences \cite{Richiardi}. 

Without any loss of generality, the GraphEM method is developed with the following setup in mind. The multi-dimensional state contains several uni-dimensional time series, each of them corresponding to a particular node. The matrix in the state equation encodes the linear dependency between the state of a given node and the previous states in the neighbor nodes. GraphEM aims at estimating the directed graph, i.e., the connections between nodes and the associated weights, that represents the (causal) dependencies between the states. GraphEM imposes a sparsity constraint in the linear matrix which promotes that each dimension of the state depends of only small subset of states in the previous time step. The method relies on the EM framework, adapted for maximum a posteriori (MAP) estimation thanks to the introduction of an $\ell_1$ penalization function. The E-step amounts to apply the Kalman filtering and RTS smoothing recursion. The M-step aims at solving a Lasso-like problem, that we address through a Douglas-Rachford proximal splitting algorithm, benefiting from great practical performance and sounded convergence guarantees \cite{CombettesDR}. Experimental results on several synthetic datasets allow to assess the fast practical convergence rate of our method, and its ability to estimate efficiently the causality graph, in particular when compared to state-of-the-art techniques relying on Granger causality.
 
The paper is structured as follows. In Section \ref{sec:format} we briefly describe the model and the filtering/smoothing algorithms, and introduce the EM framework for parameters estimation in state-space models. The novel GraphEM algorithm is presented in Section \ref{sec:proposed}. The paper concludes with a numerical example in Section \ref{sec:experiment} and some remarks in Section~\ref{sec:conclusion}.

%% file: 01_background.tex
 
\subsection{Linear state-space model}
Let us consider the following Markovian state-space model,
 \begin{equation}
   \left\{ 
   \begin{array}{ll}
       \x_{k} &= \A \x_{k-1} + \q_{k},   \\
     \y_k &= \H \x_k + \r_{k},
    \end{array}
    \right.
    \label{eq:model}
    \end{equation}
for $k=1,\ldots,K$, where $\x_{k}\in \mathbb{R}^{N_x}$ is the hidden state at time $k$,  $\y_{k}\in \mathbb{R}^{N_y}$ is the associated observation, $\A = \mathbb{R}^{N_x \times N_x}$, $\H = \mathbb{R}^{N_y \times N_x}$, $\{ \q_k \}_{k=1}^K \sim \mathcal{N}(0,\Q)$ is the i.i.d. state noise process and $\{ \r_k \}_{k=1}^K \sim \mathcal{N}(0,\R)$ is the i.i.d. observation noise process. The state process is initialized as $\x_0 \sim \mathcal{N}(\x_0 ; \bar \x_0, \P_0)$ with known $\bar \x_0$ and $\P_0$.

\subsection{Kalman filtering and smoothing}
The interest in the linear-Gaussian model of Eq. \eqref{eq:model} is usually in the computation of the sequence of filtering distributions $p(\x_k|\y_{1:k})$, where we use the shorthand $\y_{1:k} = \{ \y_j \}_{j=1}^k$. The Kalman filter provides an exact computation of those as $p(\x_k|\y_{1:k}) = \mathcal{N}(\x_k|\m_k,\P_k)$, for every $k=1,\ldots,K$~\cite{Kalman60}. The sequence of smoothing distributions (conditioned to the whole set of observations), can be also computed exactly by the RTS smoother yielding $p(\x_k|\y_{1:K}) = \mathcal{N}(\x_k|\m_k^s,\P_k^s)$. In both Kalman filter and RTS smoother, it is necessary to know the exact model, i.e., the parameters $\A$, $\H$, $\Q$, and $\R$. The explicit filter and the smoother algorithms can be found in \cite{Sarkka}.

 \subsection{EM framework for parameters estimation}

In this work, we are interested in the joint estimation of the unknown $\A$ and the hidden states. 
Parametric estimation in state-space models has been considered in the literature through three main types of methods: expectation-maximization (EM) algorithms \cite{shumway1982approach,Thiesson04}, optimization based methods \cite{Olsson}, and Monte Carlo methods \cite{Kantas09}. We will retain here the EM-based framework since it is well-suited for linear Gaussian models and it is flexible while benefiting from sounded convergence guarantees. 
EM is a method to iteratively find a maximum likelihood (ML) estimate of the parameters when the direct optimization of the posterior distribution is not feasible~\cite{MoonEM,Dempster}. It is also highly connected to majoration-minimization (MM) approaches, as it can be viewed as a special class for those~\cite{hunter2004tutorial}. The EM algorithm alternates between a majoration step consisting in building an upper bound on the neg-log-likelihood function (E-step), and the minimization of this upper bound (M-step). In the context of state-space models with parameters $\thetab$, at each iteration $i \in \mathbb{N}$ of EM method, one majorizes $- \log p(\y_{1:K} | \thetab) \leq \mathcal{Q}(\thetab;\thetab^{(i-1)})$~\cite{shumway1982approach} with $\thetab^{(i-1)}$ the parameters value of the previous EM iterate and
\begin{multline}
 \mathcal{Q}(\thetab;\thetab^{(i)}) \\
= - \int p(\x_{0:K} | \y_{1:K},\thetab^{(i)}) \log p(\x_{0:K},\y_{1:K} | \thetab) \rm{d} \x_{0:K}.
\label{eq:funQ0}
\end{multline}
Then, the M-step consists in defining $\thetab^{(i)}$ as a minimizer for $\mathcal{Q}(\cdot;\thetab^{(i)})$, so that, by construction, the sequence $(\thetab^{(i)})_{i \in \mathbb{N}}$ decreases monotonically the ML cost function, and convergence guarantees to the ML estimate can be obtained under suitable assumptions. Note that, although the EM algorithm was originally an algorithm for computing ML estimates, it can also be modified for computation of MAP estimates, as we will show in the next section.

%% file: 02_graphem.tex
 
In this section we introduce the GraphEM algorithm for the estimation of the $\A$ matrix under sparsity constraints encoded in the prior distribution $p(\A)$ (see a discussion about the prior choice in Section \ref{sec_prior}). GraphEM aims at finding the maximum of $p(\A|\y_{1:K})\propto p(\A)p(\y_{1:K}|\A)$, i.e., the MAP estimate of $\A$. 
 It is direct to show that the maximum of $p(\A|\y_{1:K})$ coincides with the minimum of $    \varphi_K(\A) = - \log p(\A) 
    - \log p(\y_{1:K}| \A)$.
In the model~\eqref{eq:model}, $\varphi_K$ can be expressed recursively for $k=1,\ldots,K$:
\begin{align}
    \varphi_k(\A) & = \varphi_{k-1}(\A) - \log p(\y_k | \y_{1:k-1}, \A)\\
    & = \varphi_{k-1}(\A) + \frac{1}{2} \log | 2 \pi \S_k(\A)| \nonumber \\ & \qquad + \frac{1}{2} \z_k(\A)^\top \S_k(\A)^{-1} \z_k(\A),
\end{align}
where $\z_k(\A) = \y_k - \H\A\m_{k-1}(\A)$  and $\S_k(\A)$ is the covariance matrix of the predictive distribution $p(\y_{k}|\y_{1:k-1})= \mathcal{N}\left(\y_{k};\H\A\m_{k-1}(\A),\S_{k}(\A)\right)$, both being side products of the Kalman filter ran for a given $\A$ (see \cite[Section 4.3]{Sarkka}). Moreover, $\varphi_0(\A) = - \log p(\A)$ is the regularization function. The direct minimization of $\varphi_K$ is made difficult due to its recursive form. The EM approach allows to construct an upper bound for it, that is more tractable and thus easier to minimize. Let us denote $\G_k = \P_k (\A')^\top (\A' \P_k (\A')^\top + \Q)^{-1}$ the output of the RTS smoother for a given $\A' \in \mathbb{R}^{N_x \times N_x}$ and set $\Sigmab = \frac{1}{K}\sum_{k=1}^K \P_k^s + \m_k^s (\m_k^s)^\top$, $\Phib  = \frac{1}{K} \sum_{k=1}^K \P_{k-1}^s + \m_{k-1}^s (\m_{k-1}^s)^\top$, $\C = \frac{1}{K} \sum_{k=1}^K \P_k^s \G_{k-1}^\top + \m_k^s (\m_{k-1}^s)^\top$. Then, following \cite[Section 12.2.3]{Sarkka}, we can show that the following function, parametrized by $\A'$, majorizes the MAP objective function $\varphi_K$ for every $\A \in \mathbb{R}^{N_x \times N_x}$:
\begin{multline}
   \mathcal{Q} (\A;\A')  =  \frac{K}{2}  \text{tr} \left(\Q^{-1} (\Sigmab - \C \A^\top - \A \C^\top + \A \Phib \A^\top) \right)  \\ + \varphi_0(\A) + \mathcal{C}
   \label{eq:majQ}
\end{multline}
where $\rm{tr}$ is the trace operator, and $\mathcal{C}$ is a constant term independent from~$\A$. The GraphEM algorithm is then constructed as described in Alg.~\ref{algo:EM_MAP}. The algorithm iterates alternating between an E-step and an M-step. In the E-step, the function $\mathcal{Q}(\A;\A^{(i-1)})$ is formed by running the Kalman filter and RTS smoother where the state matrix in the model is set to $\A^{(i-1)}$. In the M-step, the matrix $\A^{(i)}$ is updated by minimizing the latest generated majorizing function $\mathcal{Q}(\A;\A^{(i-1)})$. In the following, we provide more details about the choice of the sparsity promoting prior and the M-step implementation.

\begin{algorithm}[h!]
\noindent Initialization of $\A^{(0)}$\\
\noindent For $i = 1,2,\ldots$\\
    $\phantom{aaa}$ \textbf{(E-step)} Run the Kalman filter and RTS smoother by setting $\A':=\A^{(i-1)}$ and construct   $\mathcal{Q}(\A;\A^{(i-1)})$ using Eq. \eqref{eq:majQ}.\\    $\phantom{aaa}$ \textbf{(M-step)} Update $\A^{(i)} = \text{argmin}_A \left( \mathcal{Q}(\A;\A^{(i-1)})\right)$ (see Section \ref{sec_m_step})
        \caption{GraphEM algorithm}
    \label{algo:EM_MAP}
\end{algorithm}

\subsection{Choice of the prior}
\label{sec_prior}

Sparsity is a key feature in statistical data processing to limit the degrees of freedom in parametric models. It has been widely used for graphical model inference~\cite{Belilovsky,Friedman07}, since a sparse matrix better allows to reveal interpretable and compact network of interdependencies between the entities. Since the $\ell_0$ count measure is barely tractable, many statistical or learning approaches have invested in computable proxies, such as the $\ell_1$ norm, also called Lasso penalty. In our problem, this amounts to set $\varphi_0(\A) = \gamma \| \A\|_1$ with $\gamma$ a positive penalty parameter. The $\ell_1$ penalty is non-smooth, though convex and can be efficiently handled thanks to proximal optimization methods, as we will show hereafter. 

\subsection{Computation in the M-step}
\label{sec_m_step}
At an iteration $i \in \mathbb{N}$, the M-step minimizes the function $\mathcal{Q}(\A ;\A^{(i-1)})$, obtained by plugging $\A':=\A^{(i-1)}$ in  Eq.\eqref{eq:majQ}. This minimization problem, sometimes called Lasso regression~\cite{Tibshirani}, has been much studied in the literature of optimization \cite{Schmidt,Bach}, and most of the methods proposed to solve it rely on the proximity operator \cite{Combettes2010,Bauschke}.\footnote{See also \url{http://proximity-operator.net/}} 
 We recall that, for a proper, lower semi-continuous and convex function $f: \mathbb{R}^{N_x \times N_x} \mapsto ]-\infty,+\infty]$, the proximity operator of $f$ at $\widetilde{\A} \in \mathbb{R}^{N_x \times N_x}$ is defined as 
\begin{equation}
    \text{prox}_f(\widetilde{\A}) = \text{argmin}_\A \left( f(\A) + \frac{1}{2}\| \A - \widetilde{\A}\|^2_F \right).
\end{equation}

Let us decompose $\mathcal{Q}(\A ;\A^{(i-1)}) = f_1(\A) + f_2(\A)$, where $f_1(\A) = \frac{K}{2}  \text{tr} \left(\Q^{-1} (\Sigmab - \C \A^\top - \A \C^\top + \A \Phib \A^\top) \right)$ and $f_2 = \varphi_0$ is the prior that we have discussed in Section~\ref{sec_prior}. Function $f_1$ is quadratic, so that, for $\theta > 0$, we have:
\begin{multline*}
    \text{prox}_{\theta f_1}(\widetilde{\A}) =  \\    \text{mat} \left( 
    \left[\Id \otimes (K \Q^{-1}) + (\text{$\theta$} \Phib^{-1}) \otimes \Id \right]^{-1} \rm{vec}(K \Q^{-1}\C \Phib^{-1})
    \right),
\end{multline*}
due to the fact that the Lyapunov equation $\X \A + \A \Y = \Z$ has for solution, $\A = \rm{mat}( (\Id \otimes \X + \Y^\top \otimes \Id)^{-1} \rm{vec}(\Z))$, with $\rm{vec}$ and $\rm{mat}$ the re-ordering operations, using lexicographic order and $\otimes$ the Kronecker product. Note that if $\Q = \sigma_\Q^2 \rm{\Id}$, i.e., the state noise is isotropic, the expression simplifies into
\begin{equation*}
    \text{prox}_{\theta f_1}(\widetilde{\A}) =  \left(\frac{\theta K}{\sigma_\Q^2} \C+\widetilde{\A}\right) \left(\frac{\theta K}{\sigma_\Q^2} \Phib + \Id\right)^{-1}.
\end{equation*}
Second, the proximity operator for $f_2$ is the simple soft thresholding operator:
\begin{equation*}
    \operatorname{prox}_{\theta f_2}(\widetilde{\A}) = \left(\text{sign}(\widetilde{A}_{nm}) \times\max(0,|\widetilde{A}_{nm}| - \theta)
    \right)_{1 \leq n,m \leq N_x}
\end{equation*}

Since function $f_1$ is differentiable, a natural strategy would be to adopt the iterative thresholding method \cite{Daubechies} or an accelerated version for it \cite{Beck}. Here, following the comparative analysis from \cite{Glaudin2019}, we will prefer an algorithm that activates both terms via their proximity operator. More precisely, we will make use of the Douglas-Rachford (DR) algorithm, a fixed-point strategy for convex optimization that benefits from sounded convergence guarantees \cite{CombettesDR}, and has demonstrated its great practical performance in matrix optimization problems related to graphical inference applications~\cite{Benfenati18}. This leads to Algorithm~\ref{algo:DR}, that generates a sequence $\{\A_n\}_{n \in \mathbb{N}}$ guaranteed to converge to a minimizer of $f_1 + f_2$. In practice, for the iteration $i$ in GraphEM, we run DR method with $\theta = 1$ and initialization $\Z_0 = \A^{(i-1)}$, i.e. the majorant function tangency point. Moreover, we stop the DR loop as soon as $|(f_1 + f_2)(\A_{n+1}) - (f_1 + f_2)(\A_{n})| \leq \varepsilon$ (typically, $\varepsilon = 10^{-3}$).
\begin{algorithm}
Set $\Z_0 \in \mathbb{R}^{N_x \times N_x}$ and $\theta \in (0,2)$\\
For $n=1,2,\ldots$\\
$\phantom{aaa}$ $\A_n = \operatorname{prox}_{\theta f_2}(\Z_n)$\\
$\phantom{aaa}$ $\V_n = \operatorname{prox}_{\theta f_1}(2 \A_n - \Z_n)$\\
$\phantom{aaa}$ $\Z_{n+1} = \Z_n + \theta (\V_n - \A_n)$
\caption{Douglas-Rachford algorithm for M-step}
\label{algo:DR}
\end{algorithm}


%% file: 03_experiments.tex
The synthetic data are generated following model \eqref{eq:model} with $N_x = N_y$. A block-diagonal matrix $\A$ is considered, composed with $b$ blocks with dimensions ${b_j}_{1 \leq j \leq b}$, so that $N_y = \sum_{j=1}^b b_j$. The $j$-th diagonal block of $\A$ reads as a randomly selected auto-regressive of order one (AR1) matrix. Note that a projection on the space of matrices with spectral norm less than one is applied, so that matrix $\A$ leads to a stable Markov process. Furthermore, we set $\B = \rm{\Id}$, $K = 10^3$ and $\Q = \sigma_\Q^2 \rm{\Id}$, $\R = \sigma_\R^2 \rm{\Id}$, $\P_0 = \sigma_\P^2 \rm{\Id}$ with $(\sigma_\Q,\sigma_\R,\sigma_\P)$ some predefined values. Following this procedure, we created four datasets, whose parameters are provided in Table~\ref{tab:data}. Measures of quality to assess the performance of the methods are the relative mean square error (RMSE) on $\A$, and the precision, recall, specificity, accuracy, and F1 score for the graph edge detection, when applying a threshold of $10^{-10}$ on the entries of matrix $\A$. 

\begin{table}[h]
\centering
{\footnotesize
\begin{tabular}{|c||c|c|}
\hline
Dataset & $(b_j)_{1 \leq j \leq b}$ & $(\sigma_\Q,\sigma_\R,\sigma_\P)$\\
\hline\hline
A & $(3,3,3)$ & $(10^{-1},10^{-1},10^{-4})$\\
\hline
B & $(3,3,3)$ & $(1,1,10^{-4})$\\
\hline
C & $(3,5,5,3)$ & $(10^{-1},10^{-1},10^{-4})$\\
\hline
D & $(3,5,5,3)$ & $(1,1,10^{-4})$\\
\hline
\end{tabular}
}
\caption{Description of datasets}
\label{tab:data}
\end{table}
%
For each dataset, we ran GraphEM algorithm with the stopping rule $|\varphi_K(\A^{(i)}) - \varphi_K(\A^{(i-1)})| \leq 10^{-3}$. In practice, a maximum number of $50$ iterations is sufficient to reach this criterion. The algorithm is initialized setting $\A^{(0)}$ as an AR1 matrix. Parameter $\gamma$, balancing the weight of the sparse prior, is optimized on a single realization thanks to a manual grid search, so as to maximize the accuracy score. We also provide the results obtained by the ML estimator, also computed using EM algorithm (MLEM). In this case, the M-step has a closed form solution~\cite[Th.12.5]{Sarkka}. In addition, we compare with two Granger-causality approaches \cite{Bressler2011} for graphical modeling. The first algorithm, pairwise Granger Causality (PGC) explores the $N_x(N_x$-$1)$ possible dependencies among two nodes, at each time independently from the rest. The second approach, the conditional Granger Causality (CGC), operates similarly but for each pair of nodes, it also takes into account the information of the other $N_x-2$ signals, in order to evaluate whether one node provides information to the other when the rest of signals are observed. We also do a manual grid search for finetuning the parameters of both PGC and CGC (more information can be found in \cite{Luengo19}). Note that PGC and CGC do not estimate a weighted graph but a binary one, so that RMSE is not calculated in those cases. 

\begin{figure}
\begin{tabular}{cc}
\includegraphics[width = 0.2\textwidth]{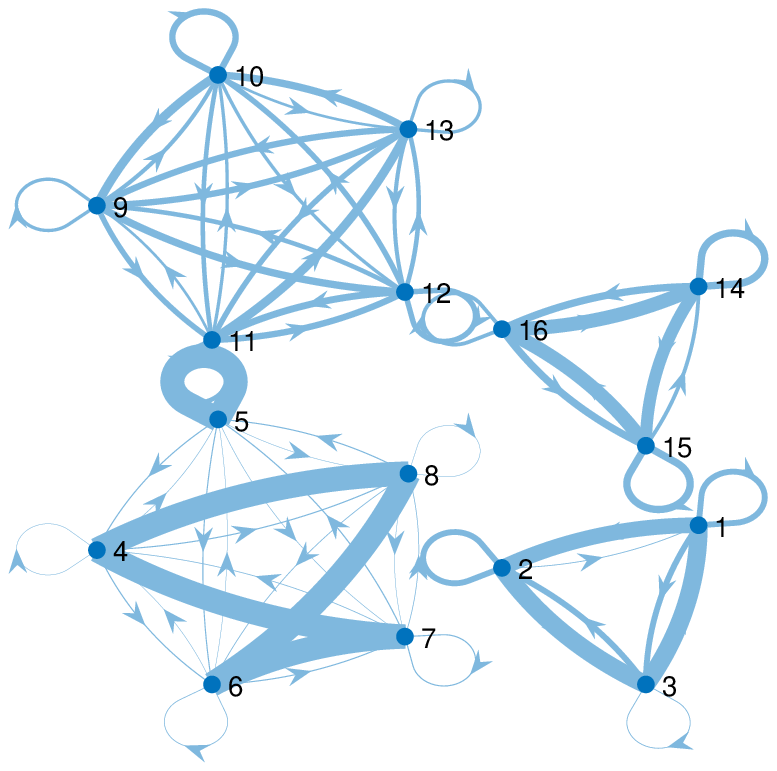} & \includegraphics[width = 0.2\textwidth]{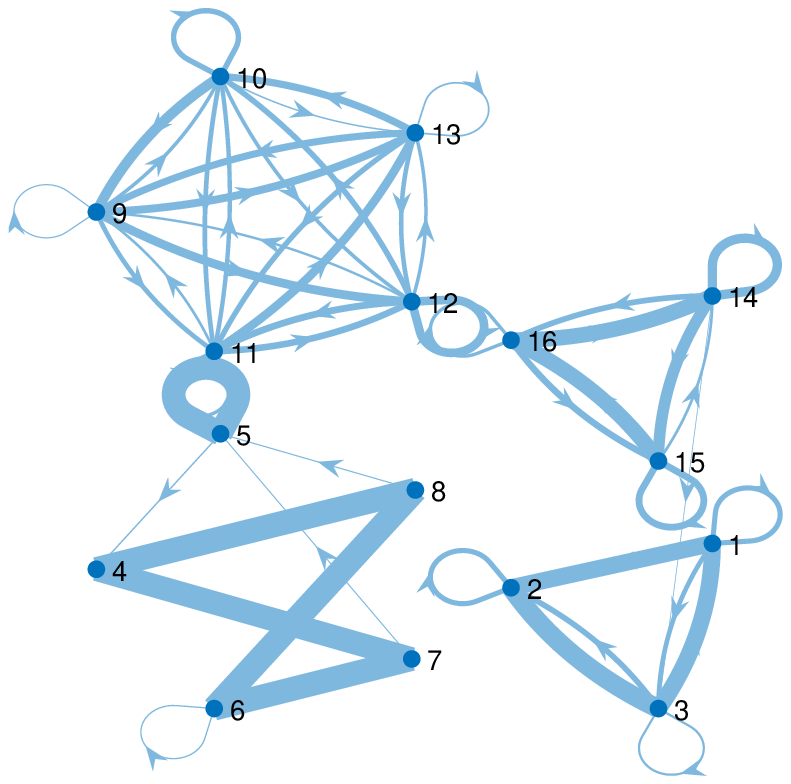} 
\end{tabular}
\vspace*{-0.2cm}
\caption{True graph (left) and GraphEM estimate (right) for dataset C.}
\label{fig:datasetC}
\end{figure}

The results, averaged on $50$ realizations, are presented in Table~\ref{tab:results}. MLEM does not promote sparsity in the graph which explains the poor results in terms of edge detectability. Moreover, GraphEM provides a better RMSE score on all examples. Regarding the graph structure, we can observe that GraphEM has also better detection scores, when compared with both PGC and CGC. We also display an example of graph reconstruction for dataset C in Fig.~\ref{fig:datasetC}, illustrating the ability of GraphEM to  recover the graph shape and weights. 


\begin{table}[t]
\centering
{\scriptsize
\begin{tabular}{|c|c||c|c|c|c|c|c|}
\cline{2-8}
\multicolumn{1}{c|}{ } & method & RMSE & accur. & prec. & recall & spec. & F1\\
\hline
\multirow{4}{*}{A} & GraphEM & $0.081$ & $0.9104$ & $0.9880$ & $0.7407$ & $0.9952$ & $0.8463$\\
  & MLEM & $0.149$ & $0.3333$ & $0.3333$ & $1$ & $0$ & $0.5$\\
  & PGC& - & $0.8765$ & $0.9474$ & $0.6667$ & $0.9815$& $0.7826$\\ 
  & CGC & - & $0.8765$ & $1$ & $0.6293$ & $1$& $0.7727$\\ 
\hline
\multirow{4}{*}{B} & GraphEM & $0.082$ & $0.9113$ & $0.9914$ & $ 0.7407$ & $0.9967$ & $0.8477$\\
& MLEM & $0.148$ & $0.3333$ & $0.3333$ & $1$ & $0$ & $0.5$\\
  & PGC & - & $0.8889$ & $1$ & $0.6667$ & $1$& $0.8$\\ 
  & CGC & - &$0.8889$ & $1$ & $0.6667$ & $1$& $0.8$\\ 
\hline
\multirow{4}{*}{C} & GraphEM & $0.120$ & $0.9231$& $0.9401$ & $0.77$ & $0.9785$ & $0.8427$\\
& MLEM & $0.238$ & $0.2656$ & $0.2656$ & $1$ & $0$ & $0.4198$\\
  & PGC& - & $0.9023$ & $0.9778$ & $0.6471$ & $0.9949$& $0.7788$\\ 
  & CGC& - & $0.8555$ & $0.9697$ & $0.4706$ & $0.9949$& $0.6337$\\ 
\hline
\multirow{4}{*}{D} & GraphEM& $0.121$ & $0.9247$ & $0.9601$ & $0.7547$ & $0.9862$ & $0.8421$\\
  & MLEM& $0.239$ & $0.2656$ & $0.2656$ & $1$ & $0$ & $0.4198$\\
  & PGC & - & $0.8906$ & $0.9$ & $0.6618$ & $0.9734$& $0.7627$\\ 
    & CGC& - & $0.8477$ & $0.9394$ & $0.4559$ & $0.9894$& $0.6139$\\ 
\hline
\end{tabular}
}
\caption{Results for GraphEM, MLEM, PGC and CGC.}
\label{tab:results}
\end{table}

